\documentclass[12pt]{amsart}
\usepackage{fancyhdr}
\usepackage[left=1.5in, right=1.5in, top=1.5in,bottom=1.5in, headsep=0.4in, footskip=0.4in]{geometry}
\usepackage{amsmath,amsfonts,amssymb}
\usepackage[T1]{fontenc}
\normalsize
\usepackage{epsfig}
\usepackage{graphicx,epsfig}
\usepackage{amsmath,amsthm,amssymb,graphicx, amscd}
\RequirePackage{amsmath}

\usepackage{amsmath,amsfonts,amsthm,amscd,amssymb,euscript,
graphics,psfrag}

\usepackage{amsfonts,amsmath,url,verbatim}

\newcommand{\afine}{{\mathbb A}}

\newcommand{\mv}{{\mathbb X}}

\newcommand{\beq}{\begin{equation}}

\newcommand{\eeq}{\end{equation}}

 \theoremstyle{definition}

\theoremstyle{plain}      
 \newtheorem{theorem} {Theorem}
\newtheorem{proposition}[theorem] {Proposition}
\newtheorem{conjecture}{Conjecture}
\newtheorem{corollary}[theorem] {Corollary}
\newtheorem{lemma} [theorem]{Lemma}

\newcommand{\Mod}{{\rm mod}}

\newcommand{\diam}{{\hbox{\rm diam}}}
\newcommand{\fp}{{\mathbb F_p}}
\newcommand{\cx}{{\mathbb C}}

\newcommand{\zed}{{\mathbb Z}}

\newcommand{\SL}{{\rm SL}}

\newcommand{\ve}{\varepsilon}

\newcommand{\ga}{\alpha}
\newcommand{\gb}{\beta}

\newcommand{\bs}{\backslash}

\newcommand{\ratls}{{\mathbb Q}}
\newcommand{\Fp}{{\mathbb F_p^*}}

\newcommand{\Fpp}{{\mathbb F_{p^2}^*}}

\newcommand{\h}{\chi}

\newcommand{\fpp}{{\mathbb F_{p^2}}}

\newcommand{\PGL}{\mbox{\it PGL}}

\newcommand{\im}{{\rm Im}}

\newcommand{\Mat}{{\rm Mat}}

\newcommand{\N}{{\rm N}}

\newcommand{\Norm}{{\rm Norm}}

\newcommand{\rot}{{\hbox{\rm rot}}}

\newcommand{\lcm}{{\hbox{\rm lcm}}}

\newcommand{\ord}{{\hbox{\rm ord}}}

\newcommand{\xx}{X^*}

\begin{document}

\title{MARKOFF SURFACES AND STRONG APPROXIMATION:  1}

\author{Jean Bourgain}
\address{IAS}
\email{bourgain@math.ias.edu}

\author{Alexander Gamburd}
\address{The Graduate Center, CUNY}
\email{agamburd@gc.cuny.edu}

\author{Peter Sarnak}
\address{IAS and Princeton University}
\email{sarnak@math.princeton.edu}

\maketitle
\section{\bf Introduction}

This is the first of three papers giving detailed proofs of the results announced in \cite{BGS16}.   The main result here is strong approximation for Markoff triples for prime moduli.  In paper two this is extended to more general moduli and applied to seiving in Markoff numbers.  In the third paper, we formulate a strong approximation conjecture for more general affine Markoff surfaces and develop the techniques to obtain similar results in these cases.

We review briefly the notation and setup from \cite{BGS16} where background and references can be found.   The Markoff surface $\mv$ is the affine surface in $\afine^3$ given by  
\beq\label{p11}
\mv  \, \, : \, \, x_1^2+x_2^2+x_3^2-3x_1 x_2 x_3=0.\eeq

The Markoff triples $\mathcal{M}$  are the positive integer solutions to \eqref{p11}.  Let $\Gamma$ be a group of affine integral morphisms of $\afine^3$ generated by the permutations of the coordiantes and the Vieta involutions $R_1, R_2,, R_3$, where $$R_3(x_1, x_2, x_3)=(x_1, x_2, 3x_1x_2-x_3)$$ and $R_2$, $R_3$ are defined  similarly.  The orbit of $(1,1,1)$ under $\Gamma$ yields all of $\mathcal{M}$ \cite{Mar79}, \cite{Mar80}.  If $\Delta$ is the group of integral morphisms generated by $\Gamma$ and the involutions which replace two of the coordinates of $x$ by their negatives, then $\mv(\zed)$ consists of two $\Delta$ orbits, namely of $(0, 0, 0)$ and $(1,1,1)$.

For $p$ a prime number the action of $\Gamma$ and $\Delta$ on $\mv$ descends to a permutation action on the finite set $\mv(\zed / p \zed)$ of solutions to \eqref{p11} in $\zed/p\zed$.  Our interest is in the orbits of this action which we often refer to as the \emph{components}.

\begin{conjecture}[Strong Approximation Conjecture] \label{Mcon}
For any prime $p$, $\mv(\zed/p\zed)$ consists of two $\Gamma$ orbits, namely $\{(0,0,0)\}$ and $\xx(\zed/p\zed) =\mv(\zed/p\zed) \bs \{(0, 0, 0)\}$
\end{conjecture}

Clearly the above conjecture implies strong approximation for $\mathcal{M}$ and $\mv(\zed)$  in the form that their reductions mod $p$ 
$$\mathcal{M} \to \xx(p)$$
and
$$\mv(\zed) \to \mv(\zed / p \zed)$$
are onto.

In particular, one of the consequences is that  the only prime congruence obstruction for a Markoff number $m$ (that is a number which is a coordinate of an $x \in \mathcal{M}$) is the one noted in \cite{Fro13}: $m\neq \frac{2}{3}, 0 \mod p$ if $p = 3 \mod 4$, $p\neq 3$.

Our first result asserts that there is a very large orbit.

\begin{theorem}\label{t1}  Fix $\varepsilon >  0$.  Then for $p$ large there is a $\Gamma$ orbit $\mathcal{C}(p)$ in $\xx(p)$ for which 
$$|\xx(p) \bs \mathcal{C}(p)|\leq p^{\varepsilon}$$
(note that $|\xx(p)| \sim  p^{2}$), 
and any $\Gamma$ orbit $\mathcal{D}(p)$ satisfies
$$|\mathcal{D}(p)| \gg (\log p)^{\frac{1}{3}}.$$
\end{theorem}

The proof of Theorem \eqref{t1} establishes the strong approximation conjecture unless $p^2-1$ is a very smooth number.   In particular, the set of primes for which the strong approximation conjecture fails is very small.

\begin{theorem} \label{t2} Let $E$ be the set of primes for which the strong approximation conjecture fails.   For $\varepsilon > 0$,  the number of primes $p \leq T$ with $p \in E$ is at most $T^{\varepsilon}$, for $T$ large.
\end{theorem}

We end the introduction by outlining the rest of the paper.  In Section \ref{s2}  we define the fundamental rotations in $\Gamma$ that are associated to an $x \in \xx(p)$ and one of its coordinates.  These act on the conic sections gotten by intersecting $\xx(p)$ with the plane corresponding to the particular coordinate.  Some basic properties of the incidence graph of the intersections of the conic sections are established.

In Section \ref{s3}, which we call the \emph{endgame},  we define the \emph{cage} $\mathcal{C}(p)$ which is shown to be a large component of $\xx(p)$.  Specifically any $x\in \xx(p)$ for which the rotation associated to one of its coordinates (see Section \ref{s2} for definitions) has order at least $p^{\frac{1}{2}+\delta_0} $ ($\delta_0>0$) is shown to be in $\mathcal{C}(p)$.

In Section \ref{s4} , \emph{the middle game}, the last statement is extended to $x$'s for which the corresponding rotation has order $p^{\varepsilon_0}$ ($\varepsilon_0 >0$ any fixed small number).  

The methods used in Sections \ref{s3} and \ref{s4} rely on  nontrivial upper bounds for the number of points lying on curves over finite fields.  In Section \ref{s3} Weil's Riemann Hypothesis \cite{Wei41} is a key tool, but this is not strong enough when the order is less than $p^{\frac{1}{2}+\delta_0}$.    In its stead we use Stepanov's auxiliary polynomial method \cite{Ste69} in the Appendix, or the recent gcd$(u-1, v-1)$ bounds of Corvaja-Zannier \cite{CZ13}, or the combinatorial method based on the projective  Szemeredi-Trotter Theorem \cite{Bou12} developed in section \ref{s4}.

Section \ref{s5}, \emph{the opening}, deals with $x$'s for which the orders of the associated rotations are  very small (for example, being uniformly bounded).
This is done by lifting the equations to characteristic zero which leads to an equation in $(\bar{\ratls})^3$ in roots of unity.  In general (that is, in the setting of Paper III) we invoke Lang's $G_m$ conjecture at this point (proven in \cite{Lau83} for example), however for the special case at hand one can show directly that $\mv^{*}(\bar{\ratls})$ has no finite $\Gamma$ orbits (the last is a necessary condition for strong approximation, in the form of Conjecture \ref{Mcon},  to hold).   This $\bar{\ratls}$ analysis leads to to part 2 of Theorem 1.

In Section \ref{s6} we assemble the various stages of our argument, explicated in Sections \ref{s3}, \ref{s4}, \ref{s5}  and give a proof of Theorem \ref{t1} in a form from which strong approximation follows if $p^2-1$ is not very smooth: see \eqref{smooth}.  To prove Theorem \ref{t2} we combine the above with a variant of the results in \cite{CKSZ14} concerning the multiplicative orders of coordinates of points of curves on $\afine^2$ over $\fp$.

\section{Preliminaries} \label{s2}

\subsection{Analysis of the conic sections}

Theorem \ref{t1}  in the weaker form that $|\mathcal C(p)|\sim |X^*(p)|$ as $p\to\infty$, can be viewed as the finite field analogue of \cite{Gol03},  where it
is shown that the action of $\Gamma$ on the compact real components of the character variety of the mapping class group of the once
punctured torus is ergodic.
As in \cite{Gol03},   our proof makes use of the rotations $\tau_{ij} \circ R_i, i\not= j$ where $\tau_{ij}$ permutes $x_i$ and $x_j$.
For example, $$\tau_{2,3} \circ R_2 (x_1, x_2, x_3) = (x_1, x_3, 3 x_1 x_3 -x_2),$$
so the action on $(x_2, x_3)$ for fixed $x_1$ is given by the rotation $\rot(3x_1)$
\beq \label{rot} \rot(3x_1)\begin{pmatrix}
x_2\\x_3\end{pmatrix} =
\begin{pmatrix}
x_3\\ 3x_1 x_3 -x_2\end{pmatrix} =
\begin{pmatrix} 0&1\\ -1 & 3x_1\end{pmatrix} \begin{pmatrix} x_2\\ x_3\end{pmatrix}.
\eeq
This rotation preserves the conic section obtained by intersecting $X^*(p)$ with the plane defined by the first coordinate being equal to the value of $x_1$;
in general, we define the conic section $C_j(a)$ as follows:
\beq \label{conic} C_j(a) = \{x_j = a\} \cap \xx(p) .\eeq

We  give an explicit description of this action,
depending on whether $x=3 x_1$ is parabolic ($x^2-4=0$, that is $x = \pm
2$), hyperbolic ($\left(\frac{x^2-4}{p}\right)=1$) or elliptic ($\left(\frac{x^2-4}{p}\right) = -1$) with $\left(\frac{\cdot}{\cdot}\right)$ being the Legendre symbol.

\begin{lemma} \label{l:par}Let $x =\pm 2$.  If
$\left(\frac{-1}{p}\right)=-1$, that is if $p \equiv -1 \mod 4$,
then $C_1(x)$ is empty.
 If
$\left(\frac{-1}{p}\right)=1$, that is if $p \equiv 1 \mod 4$,
then \beq C_{1}(2)=(2, t, t\pm 2 i), \eeq where $i^2 \equiv -1
\mod p$; \beq C_{1}(-2)=(-2, t, -t\pm 2 i), \eeq which are pairs
of disjoint lines. The action of $\rho_1=\rot(x)$ is \beq \rho_1(2,
t, t \pm 2i)=(2, t\pm 2 i, t \pm 4 i), \eeq \beq \rho_1(-2, t, -t \pm 2i)=(-2,
-t\pm 2i, -t\mp 4i), \eeq so $\rot(2)$ preserves each line and $\rot(-2)$interchanges them.
\end{lemma}

Now when $x\neq \pm 2$ we write 
$$x=\h+\chi^{-1},$$
where $\chi \in \fp$ if  $\left(\frac{x^2-4}{p}\right)=1$ and $\chi \in \fpp$ if $\left(\frac{x^2-4}{p}\right)=-1$.

Note that 
$$
\rot(x)==\begin{pmatrix} 1&1\\ \h&\frac 1 \h\end{pmatrix} \ \begin{pmatrix} \h&0\\0&\frac 1 \h\end{pmatrix}
\begin{pmatrix} \frac 1 \h& -1\\ -\h&1\end{pmatrix} 
\Big(\frac 1 \h -\h\Big)^{-1} =\begin{pmatrix} 1&1\\ \h&\frac 1 \h\end{pmatrix}
\begin{pmatrix} \h&0\\ 0&\frac 1 \h\end{pmatrix} \begin{pmatrix} 1&1\\ \h&\frac 1 \h\end{pmatrix}^{-1}
$$
and consequently 
$$
\rot(x)^\ell =\begin{pmatrix} 1&1\\ \h&\frac 1 \h\end{pmatrix} \begin{pmatrix} \h^\ell & 0\\ 0&\frac 1{\h^\ell}\end{pmatrix}
\begin{pmatrix} 1&1\\ \h&\frac 1 \h\end{pmatrix}^{-1} =\begin{pmatrix} 1&1\\ \h&\frac 1 \h\end{pmatrix}
\begin{pmatrix} \h^\ell &0\\ 0&\frac 1{\h^\ell}\end{pmatrix} \begin{pmatrix}\frac 1 \h&-1\\ -\h&1\end{pmatrix}\Big( \frac 1 \h-\h\Big)^{-1}
$$
and 
$$
\langle \rot(x) \rangle =\Big (\frac 1 \h-\h\Big)^{-1} \left\{ \begin{pmatrix} \frac {\h_1} \h-\frac \h{\h_1}& \frac 1{\h_1}-\h_1\\
\h_1- \frac 1{\h_1}  & \frac 1{\h_1}-\h \h_1\end{pmatrix} ; \h_1\in \langle \h\rangle\right\}.
$$
Consequently $C_1(x)$ contains all elements

$$
\left(\Big(\h-\frac{1}{\h}\Big)^{-1} \left( 
\Big(x_3-\frac{x_2}{\h}\Big)\h_1+ (\h x_2 -x_3)\frac{1}{\h_1}\right),
\Big(\h-\frac{1}{\h}\Big)^{-1} \Big((\h x_3-x_2)\h_1+\Big(x_2-\frac{x_3}{\h}\Big) \frac{1}{\h_1}\Big)\right)
$$

with $\h_1\in\langle \h\rangle$.

Note that 
$$
\text {Proj}_{x_2}(C_x)\supset \Big\{a \h_1+ \frac{b}{\h_1}; \h_1\in\langle \h\rangle \Big\}
$$
where
$$
a=\Big(\h-\frac{1}{\h}\Big)^{-1} \Big(x_3 -\frac{x_2}{\h}\Big) \ \text { and } \ b=\Big(\h-\frac{1}{\h}\Big)^{-1} (\h x_2-x_3)
$$
satisfy
$$
\sigma = ab=\Big(\frac x{\h-\frac 1\h}\Big)^2 =\Big( \frac {\h+\frac 1 \h}{\h-\frac 1\h}\Big)^2 \not= 1.
$$

Denoting by $\rho$ a primitive root of $\fp$, a hyperbolic element
$x$ can be written in the form \beq x= \rho^j +\rho ^{-j}. \eeq
For a hyperbolic element we  let $\ord(x) =\frac{p-1}{j}$.

 An
elliptic element $x$ can be written in the form \beq x= \xi^{j} +
\xi^{-j}, \eeq where $\xi$ is an element in $\fpp$, 
$\xi=(\tilde{\rho})^{p+1}$, where $\tilde{\rho}$ is a generator or
the multiplicative group of $\fpp$.  For an elliptic element we let
$\ord(x) =\frac{p+1}{j}$.

\begin{lemma} \label{l:hyp} Let $x$ be hyperbolic; write \beq x= w+w^{-1},\eeq
where $w=\rho^j \in \fp$.  Then $C_1(x)$ is a hyperbola with $p-1$
points. Set \beq \kappa(x)=\frac{x^2}{x^2-4}.\eeq Let \beq H(x) =
\{\left(t, \frac{\kappa(x)}{t}\right ) \, | \, t\in \fp^{*} \}.
\eeq Then we have the following map from $H(x)$ to $C_1(x)$: \beq
\left(t, \frac{\kappa(x)}{t}\right) \to \left(x,
t+\frac{\kappa(x)}{t},
 tw+\frac{\kappa(x)}{tw}\right). \eeq
 In these coordinates
 \beq \rho_1\left(t, \frac{\kappa(x)}{t} \right)=\left(tw, \frac{\kappa(x)}{tw}
 \right). \eeq
\end{lemma}

\begin{lemma} \label{l:ell} Let $x$ be elliptic; write \beq x= v+v^{p},\eeq
where $v\in \fpp - \fp$, $v^{p+1}=1$.    Then $C_1(x)$ is an
ellipse with $p+1$ points. Set \beq
\kappa(x)=\frac{x^2}{x^2-4}.\eeq Let \beq E(x) = \{(t, t^p ) \, |
\, t\in \fpp, t^{p+1}=\kappa(x) \}. \eeq Then we have the following
map from $E(x)$ to $C_1(x)$: \beq (t, t^p) \to \left(x,
t+\frac{\kappa(x)}{t},  tv+\frac{\kappa(x)}{tv}\right). \eeq
 In these coordinates
 \beq \rho_1 \left(t,
\frac{\kappa(x)}{t}\right) =\left(tv, \frac{\kappa(x)}{tv}
 \right).\eeq
\end{lemma}

\subsection{Incidence graph for the conic sections}  We treat the case of $p\equiv 3  (\Mod 4)$ (the case of $p\equiv 1  (\Mod 4)$ is simpler  because of the special point in Lemma \ref{l:par}).  Let $\xx(p)$ be the Markoff triples $\Mod \, p$; $\xi$ any coordinate of a triple, $\xi \neq 0, \pm \frac{2}{3}$.

For $j\neq k$ and $\xi, \eta$
\beq \label{p21} 
|C_j(\xi)\cap C_k(\eta)|=0, 1, 2. \eeq
To determine which  it is, the intersection consists of all $z$'s such that 
\beq \label{p22}
\xi^2+\eta^2+z^2= 3 \xi \eta z,\eeq
which  has a solution if
$$9 \xi^2 \eta^2 - 4(\xi^2 + \eta^2)$$
is a square in $\fp$.

In terms of Legendre's symbol \beq \label{p23}
|C_j(\xi) \cap C_k(\eta)|= 1+\left(\frac{9 \xi^2 \eta^2 - 4(\xi^2 + \eta^2)}{p}\right). \eeq

So each $C_j(\xi)$ meets $\frac{p-1}{2}$ $C_k(\eta)$'s.   Define \emph{the incidence graph} $I(p)$  of $\xx(p)$ to have vertices $C_j(\xi)$'s  with the number of edges between $C_j(\xi)$  and  $C_k(\eta)$  being $|C_j(\xi) \cap C_k(\eta)|$.

\begin{proposition} \label{p2} For $p$ large ($p>10$) the incidence graph $I(p)$ is connected and in fact $\diam(I(p))=2$.
\end{proposition}

\proof   Fix $\xi_1, \xi_2$ and $i, j$, say $i, j \in \{1,2\}$.  We seek $y \in \fp$ such that $C_3(y) \cap C_i(\xi_1)\neq \phi$ and 
$C_3(y) \cap C_j(\xi_2)\neq \phi$.  This leads to solving the pair of equations:

\beq \label{p24} 
\left\{
\begin{aligned}
&(9 \xi_1^2 -4) y^2-\lambda^2 = 4 \xi_1^2\\
&(9 \xi_2^2 -4) y^2-\mu^2 = 4 \xi_2^2
\end{aligned}
\right.
\eeq
for $y, \lambda, \mu \in \fp$.
In $\xi_1^2= \xi_2^2$ then \eqref{p24} reduces to the first equation (take $\lambda=\mu$) and since $9 \xi_1^2- 4 \neq 0$ and $\xi_1^2\neq 0$, it defines a conic section.  Thus for $p$ large  it has a solution and provides us with our $y$.   If $\xi_1^2 \neq \xi_2^2$ then \eqref{p24} defines an irreducible curve in $\afine^3$.
  Thus again for $p$ large it has solutions over $\fp$ providing us with our desired $y$.  It follows that the distance in $I(p)$ between any two points is at most $2$.   On the other hand,   $C_i(\xi)$ and $C_i(\eta)$ are not joined if $\xi \neq \eta$ nor is $C_i(\xi)$ joined to half of the $C_j(\eta)$'s for $j\neq i$.   Hence $\diam( I(p)) =2$. 
\qed

\section{Endgame} \label{s3}

\subsection{Use of Weil's bound}  
We begin with the following 
\begin{proposition} \label{p31}
If $x=(x_1, x_2, x_3)$ is in $\xx(p)$ and for some $j \in \{1, 2, 3\}$ the order of the induced rotation $\rot(x_j)$ is at least $p^{\frac{1}{2}+\delta}$ ($\delta>0$ fixed) then $x$ is joined to a point $y$ in $\xx(p)$ one of whose induced rotations is of maximal order.\end{proposition}

\proof Consider first the case that $x_1$ (say $j=1$) is hyperbolic.   In light of the discussion in section \ref{s2}, $x=(x_1, x_2, x_3)$ is connected to the points in $\xx(p)$ of the form 
\beq \label{27}
(x_1, \ga_1 t +\ga_2t^{-1}, \ga_3 t +\ga_4 t^{-1}) \eeq
with $t\in H$, a cyclic subgroup of $\fp^{*}$. Here $|H| \mid (p-1)$; we set $e_H=\frac{p-1}{|H|}$.
Our aim is to produce $t$'s in $H$ for which there is a primitive root $y\in \fp^{*}$ satisfying
\beq \label{28}
\ga_1 t +\ga_2t^{-1} = y+y^{-1}.\eeq
Let $P(H) (=P_{\ga_1, \ga_2}(H))$ denote the number of such solutions.

A subgroup $K$ of $\fp^{*}$ is determined by its order $|K|$ which divides $p-1$; let $d_K=\frac{p-1}{|K|}$.  Let $f_{H}(K)=f_{H}(d_K)$ be the number of solutions to 
\beq \label{p29} \ga_1 t +\ga_2t^{-1} = y+y^{-1} \, \, , t\in H\, \, , y\in K\eeq
(note that the traces of the matrices 
that we produce, namely the common values of the
left- and right-hand side of \eqref{p29},  are hit with multiplicity $2$
in both $t$  and $y$  in our counting).
Clearly
\beq \label{p30}
f_{H}(K) \leq 2 \min(|K|, |H|). \eeq
We can estimate $f_{H}(K)$, at least if $|H|\geq p^{\frac{1}{2} +\delta}$, using Weil's Riemann Hypothesis for curves over finite fields.  The map $$\xi \to \xi^{d_{K}}, \, \, \eta \to \eta^{e_H}$$
sends solutions of \beq \label{31} 
C_{\ga_1, \ga_2} \, :\, \ga_1 \eta^{e_H} +\ga_2 \eta^{-e_H} =\xi^{d_K}+\xi^{-d_K} \eeq
to solutions of \eqref{p29} and it is $e_H d_K$ to $1$.  Hence if $N(C_{\ga_1, \ga_2})$ is the number of solutions to \eqref{31} then 
\beq \label{32} f_H(K) = \frac{N(C_{\ga_1, \ga_2})}{e_H d_K}. \eeq
  As we prove below (see Lemma \ref{irred}) the curve  $C_{\ga_1, \ga_2}$ i is absolutely irreducible.   Since its genus  is $O(e_H d_K)$,  applying Weil bound  yields \beq \label{33} N(C_{\ga_1, \ga_2})=p+ O(\sqrt{p} e_H d_K). \eeq
Hence
\beq \label{34}
f_{H}(K) = \frac{p}{e_H d_K} + O(\sqrt{p}). \eeq
By inclusion/exclusion
\beq \label{35} P(H)=\sum_{d \mid (p-1)} \mu(d) f_{H}(K_d), \eeq
where $\mu$ is the Mobius function.
Hence 
\beq \label{36} P(H)=\sum_{d|(p-1)} \mu(d) \left( \frac{|H|}{d} + O(\sqrt{p})\right) = |H|\sum_{d|(p-1)} \frac{\mu(d)}{d} + O_{\varepsilon} (p^{\frac{1}{2} + \varepsilon}) = |H|\frac{\varphi(p-1)}{p-1} + O_{\varepsilon} (p^{\frac{1}{2} + \varepsilon}).\eeq
Here $\varphi$ is the Euler function and it satisfies $\varphi(n) \gg_{\varepsilon} n^{1-\varepsilon}$ and hence from \eqref{36} we deduce that $P(H)>1$ under the assumption that $|H| \geq p^{\frac{1}{2} +\delta}$.  This proves Proposition \ref{p31} in the hyperbolic case.

Now consider the case of $x$ elliptic.  Let $D$ be a non-square element in $\fp$.   Then $\fpp = \fp[\sqrt{D}]$  and we can parametrize the subgroup $H_1$ as follows
\beq \{(\xi+\sqrt{D} \eta)^{d_1}; \xi, \eta \in \fp; \xi^2-D\eta^2 =1\},\eeq
where $d_1=\frac{p+1}{|H_1|}$.  The conic section $C_1(x)$ is an ellipse which can be parametrized as 
\beq \alpha^2-D \beta^2= \kappa(x),\eeq
where $\kappa(x)=\frac{x^2}{x^2-4}.$   We seek $\alpha$ which can be written as $\alpha=u+u^{-1}$ with $u$ a primitive root in $\fp^{*}$.

Now
\beq (\xi +\sqrt{D})^n = g_n(\xi) + h_n(\xi) \sqrt{D}, \eeq
where $g_n$, $h_n$ are integral polynomials 
\beq g_n(\xi) =\sum_{i=0}^{\lfloor n/2\rfloor} \binom{n}{2i} D^i \xi^{n-2i},\eeq
\beq h_n(\xi) =\sum_{i=0}^{\lfloor n/2\rfloor} \binom{n}{2i+1} D^i \xi^{n-2i+1}.\eeq

Let \beq g_n(\xi, \eta ) =\sum_{i=0}^{\lfloor n/2\rfloor} \binom{n}{2i} D^i \xi^{n-2i} \eta^{2i},\eeq
\beq h_n(\xi, \eta ) =\sum_{i=0}^{\lfloor n/2\rfloor} \binom{n}{2i+1} D^i \xi^{n-2i+1} \eta^{2i-1} .\eeq

Then we have 

\beq (\xi +\sqrt{D} \eta)^n = g_n(\xi, \eta) + h_n(\xi, \eta)\sqrt{D} .\eeq

Now we seek to bound $f(H_1, K)$   with $K$ subgroup of   $\fp^{*}$ with $d_2=\frac{p-1}{|K|}$.   As in the hyperbolic case this is  given by $\frac{M(d_1, d_2)}{d_1 d_2}$ where $M(d_1, d_2)$ now counts the number of points on the following curve in $\fp^3$:

\beq \label{nonsplit}
\left\{
\begin{aligned}  &\xi^2 -D \eta^2 =1 \\&  g_{d_1}(\xi, \eta) = \mu^{d_2} +\mu^{-d_2} \end{aligned} \right. \eeq

This is a curve of genus  $O(d_1 d_2)$ and we can apply Weil bound and inclusion-exclusion as in the hyperbolic case to produce the primitive $u$.  This completes the proof of Proposition \ref{p31}.
\qed

\begin{lemma} \label{irred}   Suppose $\ga_1 \ga_2 \neq 1 \mod p$ .  Then the curve 
$$\ga_1 y^{e} + \ga_2 y^{-e} = x^d + x^{-d}$$ 
is absolutely irreducible.
\end{lemma}

\proof  Consider 
$$P(x,y) = \ga_1 x^d y^{2e} +\ga_2 x^d - x^{2d} -y^e \in \overline{\fp}[X,Y].$$   
For $d=1$  $P$ is clearly irreducible.  Let $d>1$ and assume that $P$ is not irreducible and $f(x,y)=\sum a_{jk} x^j y^k \in \overline{\fp}[X,Y]$ an irreducible factor.
Assume $d\geq e$ and $u$ a $d$-th root of unity.  Since for $0\leq s \le d$, $P(x,y) = P(u^s x,y)$, also 
$$f_s(x,y)=f(u^s x,y) =\sum a_{jk} u^{sj} x^j y^k$$
is an irreducible component of $P$.  Thus either  $f_s$ and $f_{s'}$ are coprime or $f_s \sim f_{s'}$.  
Since $a_{0k} \neq 0$ for some $k$ (otherwise $x$ would be a factor of $f(x,y)$), it follows that $f_s=f_{s'}$ if $f_s \sim f_{s'}$.

\textbf{Case 1} The $f_s$ are not pairwise coprime.

The $f(x,y) = f_s(x,y)$ for some $0< s<d$, implying that $u^{s_j}=1$, i.e. $s_j \equiv 0 (\mod d)$ if $a_{jk} \neq 0$.  It follows that $d$ has a divisor $d_1>1$ such that $d_1 | j$ if $a_{jk} \neq 0$ and hence $f(x,y)$ has form $f(x,y) = g(x^{d_1}, y)$.   The polynomial $g(x,y)$ is therefore a factor of 
$Q(x,y)= \ga_1 x^{d_2} y^{2e} + \ga_2 x^{d_2} - x^{2 d_2} y^{e} - y^e$ with $d_2= \frac{d}{d_1}$ and we lowered $d$.

\textbf{Case 2} The $f_s$ $(0 \leq s \le d)$ are mutualy coprime.
Define $$P_1(x,y)= \prod_{s=0}^{d-1} f_s(x,y),$$
which divides $P$.  Degree considerations show that 
$$2d \geq d \deg_x  f, \, \, 2e \geq d \deg_y f, \, \, 2d+e \geq d  \deg f.$$

\textbf{Case 2.1}  $\deg_x f >1, \, \, \deg_y f> 1.$

It follows that $\deg_x f = 2$, $\deg_y f = 2$, $e=d$, $\deg f=3$, and $P(x,y) = P_1(x,y)$.  With $u$ as above, $\varphi(x,y) = f(x, uy)$ is an irreducible factor of $P(x,y)$.  Therefore for some $0\leq s \le  d$
$$\varphi(x,y) = \sum_{j, k \leq 2} a_{jk} u^k x^j y^k \sim f_s(x,y) = \sum a_{jk} u^{sj} x^j y^k.$$
Consequently, there is some $0\leq l \le d$ such that $k-sj\equiv l (\mod d)$ if $a_{jk} \neq 0$.  Since
$$\ga_1 x^d y^{2d} + \ga_2 x^d - x^{2d}y^d - y^d = \prod_{0 \leq s \le d} f_s(x,y),$$
clearly $a_{0,1} \neq 0$, $a_{1,0} \neq 0$, and therefore $1\equiv l \equiv -s (\mod d)$, i.e. $k+j = 1 (\mod d)$ if $a_{jk} \neq 0$.   Since $\deg f =3$, $2 \equiv 0 (\mod d)$, hence $d=2$ and $a_{1,1} =a_{2,0} = a_{0,2} =0$.
Thus
\begin{equation*}
\begin{aligned}
\ga_1 x^2 y^4 +
ga_1 x^2 - x^4 y^2 -y^2 &\sim(a_{21} x^2 y +a_{12} x y^2 + a_{10} x + a_{01} y) (a_{21} x^2 y - a_{12} xy^2 - a_{10}x + a_{01}y)\\
&\sim y^2(a_{21} x^2 +a_{01})^2 - x^2 (a_{12} y^2 +a_{10})^2.
\end{aligned}
\end{equation*}

Setting $a_{0,1}=1$ gives 
$-y^2(a_{21}x^2+1)^2 + x^2(a_{12} y^2 + a_{10})^2$ and $a_{21}^2=1$, $a_{21}-a_{12}a_{10} =0$, $a_{12}^2 =\ga_1$, $a_{10}^2=\ga_2$.
But this contradicts the assumption $\ga_1 \ga_2 \neq 1$.

\textbf{Case 2.2}  $\deg_x f=1$ or $\deg_y f =1$.

Assume $\deg_y f =1$, say.  Then there are coprime $a(x), b(x) \in \overline{\fp}[X]$ such that $P(x, \frac{a(x)}{b(x)})=0$, that is 
$$\ga_1 x^d a(x)^{2e} + \ga_2 x^d b(x)^{2e} - x^{2d} a(x)^e b(x)^e - a(x)^eb(x)^e= 0.$$ Since $a(x), b(x)$ are coprime, it follows that  $a(x)^e | x^d$, $b(x)^e | x^d$, hence $a(x)$ or $b(x)$ is constant.  If, say, $b(x)$ is constant, previous equation implies $x^d | a(x)^e$, hence $a(x)^e = \gamma x^d$ and 
$$\ga_1  \gamma^2 x^{2d} + \ga_2 b^{2e} - \gamma b^2 x^{2d} - \gamma b^2 =0.$$
It follows that $\ga_1 \gamma = b^2$, $\ga_2  b^2 =\gamma$, hence $\ga_1 \ga_2 =1$, contradicting the assumptions that $\ga_1 \ga_2 \neq 1$.
\qed

\subsection{The Cage}
A point $x=(x_1, x_2, x_3) \in \xx(p)$ is called maximal if $\ord(\rot ( x_j))$ is maximal for some $j$.  Note that the condition that the order of $\rot(3 x_j)$ be maximal depends only on $x_j$  and not on the other coordinates of $x$ (since it depends on the order of $\lambda$, where $\lambda + \lambda^{-1} = 3 x_1$ in $\Fp$ of $\Fpp$).   We call $\xi \in \fp$ maximal if it is of maximal order.

By the \emph{cage} we mean a set of maximal elements in $\xx(p)$.  We claim that the cage is connected, that is to say if $\hat{x}$ and $\hat{y}$ are in the cage then $\hat{x}$ is connected to $\hat{y}$.   Let $\xi$ be the coordinate of maximal order of $\hat{x}$ and $\eta$ be the coordinate of maximal order of $\hat{y}$, so that $\hat{x}$ is connected to all points in
 $C_j(\xi)$ and similarly $\hat{y}$ is connected  to all the points in $C_k(\eta)$.

Now according to Proposition \ref{p2}, which when extended with and inclusion/exclusion argument as in Proposition \ref{p31} gives a $y$ of maximal order such that $$P \in C_j(\xi)\cap C_l(y);$$
$$Q \in C_k(\eta) \cap C_l(y).$$
Since $y$ is maximal,  $P$ and $Q$ are connected by $\Gamma$ and $P$ is $\Gamma$- connected to $\hat{x}$ and $Q$ to $\hat{y}$.   Thus $\hat{x}$ is $\Gamma$ connected to $\hat{y}$.

Denote by $\mathcal{C}(p)$ the connected component of $\xx(p)$ (under the $\Gamma$ action) that contains the cage, then $\mathcal{C}(p)$ is our large component.

\section{Middle Game} \label{s4}  In the endgame (section \ref{s3}) we connected any $x \in \xx(p)$ of order  (that is $\max (\ord (\rot(x_j)))$) $l \geq p^{\frac{1}{2}+\delta_0}$ ($\delta_0>0$) to the cage in one step.  In this section we allow any number of moves to do the connecting.  In particular  any $x$ of order at least $p^{\varepsilon}$ is shown to be in the giant component.  As in section \ref{s3} the $y$'s which are joined to a given $x$ whose order is $l_1$ via the corresponding rotation and which have orders $l_2$ (here $l_1$ and $l_2$ divide $p-1$ or
$p+1$) correspond to solutions of an equation (with $\sigma \in \fp$, $\sigma \neq 1$):
\beq \label{p41}
\left.
\begin{aligned}
&h_1+\frac{\sigma}{h_1}=h_2+\frac{1}{h_2}, \sigma \not= 1\\
&\text {with $h_1 \in H_1,  h_2  \in H_2$ with $H_1, H_2$ subgroups of $\fp^{*}$ or $\fpp^{*}$}.
\end{aligned}
\right\}
\eeq

In \eqref{p41} we have $|H_1|=l_1, \, |H_2|=l_2$.


If we have an upper bound on the number of solutions to \eqref{p41} so that on summing over all $l_2 \leq l_1$ with $l_2$ dividing $p-1$ or $p+1$, yields a quantitiy which is less than $l_1$, then there is at least one $h_1$ for which the corresponding $y$ will have order bigger than $l_1$.   We then repeat this procedure replacing $x$ by this $y$ and so on, until the order of the element is at least $p^{\frac{1}{2}+\delta_0}$.  At that point we are in the endgame and can finish.

The key therefore is a suitable upper bound to the number of solutions to  \eqref{p41}.  
Our original treatment used Stepanov's technique of auxiliary polynomials in his elementary proof of this Riemann Hypothesis for curves.  This yields explicit and reasonably sharp estimates which  are ample for our application.  We carry this out in the   Appendix partly to illustrate the flexibility of this method.

Subsequently the recent powerful technique for estimating from above the g.c.d. of $(u-1, v-1)$, of Corvaja and Zannier \cite{CZ13} yields sharper bounds.  This is relevant for the purpose of giving efffective bounds on $p$ after which Theorem \ref{t1} takes effect.  The precise upper bound to \eqref{p41} established by \cite{CZ13} is 
$$
20\max\Big\{(|H_1|.|H_2|)^{1/3}, \frac {|H_1|.|H_2|}p\Big\}.
$$

The third treatment,  and the one which we develop in this section,  while  special to \eqref{p41},  is robust in that the upper bound requires little further structure and it is  suitable for generalisation for more general moduli as demonstrated in Paper II. It is  based on the following projective Szemeredi-Trotter theorem  (Proposition 2 in \cite{Bou12}).

\medskip

\begin{theorem} \label{sth}
Let $\Phi:\mathbb F_p\to \text {\, \Mat}_2 (\mathbb F_p)$ such that $\det \Phi$ does not vanish identically and $\im \Phi\cap \PGL_2(\fp)$ is not contained in 
a set of the form $\mathbb F_p^* \cdot gH$ for some $g\in \SL_2(\fp)$ and $H$ a proper subgroup of $\SL_2(\fp)$. Then the following holds.

Given $\ve>0, r>1$, there is $\delta>0$ such that if $A\subset P^1(\fp)$ and $L\subset \mathbb F_p$ satisfy
\beq \label{26} 
 1\ll |A|< p^{1-\ve}\eeq
\beq \label{27} 
 \log |A|< r\log |L|.\eeq
Then
\beq \label{28}
|\{(x, y, t)\in A\times A\times L; y=\tau_{\Phi(t)} (x)\}|< |A|^{1-\delta} |L|,
\eeq
where for $g=\begin{pmatrix} a&b\\ c&d\end{pmatrix}$, $\tau_g(x)=\frac{ax+b}{cx+d}$.
\end{theorem}

Using Theorem \ref{sth} we prove the following:

\begin{proposition} \label{key} Given $\delta>0$ there is $\tau<1$ and $C_{\tau}$ depending on $\delta$ such that if $p^{\delta} < |H_1|< p^{1-\delta}$ then the number of solutions to 
\eqref{p41} is at most $C_{\tau} |H_1|^{\tau}$.
\end{proposition}

\proof For $h \in H$, a subgroup of  $\mathbb F_p^*$ or $\mathbb F_{p^2}^*$,  denote by 
$\tilde{h} =h+h^{-1}$.  Similarly we denote by $\widetilde{H}=\{\tilde{h} \, | \, h\in H\}$.

Suppose that \eqref{p41} has $T$ solutions.   Then 

\beq \label{p41a}
\left.
\begin{aligned}
&h_1+\frac{\sigma}{h_1}=u\\
&h_1 t +\frac{\sigma}{h_1 t}=v,
\end{aligned}
\right\}
\eeq
where $h_1, t \in H_1$ and $u, v \in \widetilde{H_2}$, has at least $T^2$ solutions. 

 Elimination of $h_1$ in \eqref{p41a}  yields

\beq \label{22}
u^2+v^2 -\Big( t+\frac 1t\Big) uv +\sigma\Big(t-\frac 1t\Big)^2 =0\eeq
which, by assumption,  has at least $T^2$ solutions in $(t, u, v)\in H_1\times \widetilde{H_2}
\times \widetilde{H_2}$.

Next, let $u=\tilde{f_1}$, $v=\tilde{f_2}$ with $f_1$ and $f_2\in H_2$ and define the following elements $x, y \in \widetilde{H_2}$:
$$
x= \widetilde{(f_1 f_2)} = f_1 f_2 +\frac{1}{f_1 f_2},
$$
$$
 y= \widetilde{(f_1 f_2^{-1})}=\frac{f_1}{f_2}+\frac{f_2}{f_1}.
$$
Thus $uv=x+y, u^2+v^2=xy +4$ and equation \eqref{22}  gets transformed into

\beq \label{23} 
xy- \tilde t (x+y)+ \sigma (\tilde t)^2 +4(1-\sigma)=0.\eeq
Denoting
\beq \label{24} 
\alpha =\tilde t \ \text { and } \ \beta =\sigma(\tilde t)^2 +4(1-\sigma)\eeq
we obtain
\beq \label{25} 
y =\frac {\alpha x-\beta}{x-\alpha} =\tau_g(x)\eeq
with
$$
g=\begin{pmatrix} \alpha&\beta\\ 1 & -\alpha\end{pmatrix} = g(\tilde t)
$$
and $\tau_g$ the Mobius transformation.

Equation \eqref{25}  has at least $T^{2}$ solutions in $(x, y, \tilde t) \in 
\widetilde H_2\times \widetilde H_2 \times \widetilde H_1$.

We apply Theorem \ref{sth} taking $\Phi (t)$
\beq \label{29}
\Phi(t) =
\begin{pmatrix}
t & -\sigma t^2-4(1-\sigma)\\ 1& -t\end{pmatrix}
\eeq
and $A=\widetilde H_2, L=\widetilde H_1$.  

 We verify the assumption on $\Phi$.
Since $\sigma\not= 1$, det\,$\Phi(t) =(1-\sigma)(4-t^2)$ does not vanish identically.
It remains to show that
$$
\Big\{ \Phi(s)^{-1}\Phi(t)\ \frac {\det \Phi(s)}{\det\Phi(t)}; s, t \in\mathbb F_p\}
$$
is not contained in a proper subgroup $H$ of $\SL_2(p)$.

By \eqref{29}
$$
\begin{aligned}
\Phi (s)^{-1} \Phi(t) &=\begin{pmatrix}
-s& \sigma s^2+4(1-\sigma)\\ -1& s\end{pmatrix} \begin{pmatrix} t&-\sigma t^2-4(1-\sigma)\\ 1& -t\end{pmatrix}
\\[12pt]
&= \begin{pmatrix}
-st+\sigma s^2 +4(1-\sigma)\  & \ (s-t)(4(1-\sigma)-\sigma st)\\
s-t & -st+\sigma t^2 +4(1-\sigma))\end{pmatrix}
\end{aligned}
$$
Taking
$$
s=\sigma t + \frac {4(1-\sigma)}t
$$
gives
\beq \label{210}
(1-\sigma) \Big(\frac 4t-t\Big)
\begin{pmatrix} \sigma (1+\sigma)t+\frac {4(1-\sigma)} t \ &  \ - \sigma^2t^2 +4(1-\sigma)^2\\
1& 0\end{pmatrix} = (1-\sigma) \Big(\frac 4t -t\Big) g_t.\eeq

As the proper subgroups of $\SL(2, \fp)$ have trivial second commutator \cite{Suz82}, it suffices  to show that
\beq \label{211}
(g_{t_1} g_{t_2} g^{-1}_{t_1} g^{-1}_{t_2})(g_{t_3} g_{t_4} g_{t_3}^{-1} g_{t_4}^{-1})(g_{t_2} g_{t_1} g_{t_2}^{-1} g_{t_1}^{-1})
(g_{t_4} g_{t_3} g_{t_4}^{-1} g_{t_3}^{-1})\eeq
is not identically one for $t_1, t_2, t_3, t_4 \in \mathbb F_p^*$.
If this were the case, the same would be true for $t_1, t_2, t_3, t_4$ taken in an extension field of $\mathbb F_p$ so as to make
\beq \label{212}
t^2 =\frac {4(1-\sigma)^2-\ve}{\sigma^2}  \ (\ve =\pm 1)\eeq
solvable.

Taking $t=\pm \kappa$ satisfying \eqref{212}, we get
\beq \label{213}
g_{\pm\kappa}=\begin{pmatrix} \pm \frac \sigma\kappa [(1+\sigma) \frac {4(1-\sigma)^2-\ve}{\sigma^2} +4(1-\sigma)]& \ve\\
1& 0\end{pmatrix}.\eeq
We choose $\ve=\pm 1$ as to ensure that 
$$
(1+\sigma)(4(1-\sigma)^2-\ve)+4\sigma^2(1-\sigma)\not= 0
$$ 
and obtain matrices
\beq \label{214}
g_{\pm} =\begin{pmatrix} \pm\eta&\ve\\ 1&0\end{pmatrix}\eeq
that clearly generate $\SL_2(p)$. 

Consequently, Theorem \ref{sth} is applicable, yielding the bound $T^2  \ll |H_2|^{1-\tau} |H_1|$.

\qed

\section{Opening} \label{s5}  
The analysis of the previous sections shows that we can connect $x \in \xx(p)$ whose order is at least $p^{\varepsilon}$ (or smaller if the divisors of $p^2-1$ are not too numerous) to the cage.  To deal with all $x$'s and in particular ones whose orders are uniformly bounded (indepndent of $p$) we lift to characteristic zero.  In this connection we observe first that if the action of $\Gamma$ on $\xx(\bar{\ratls})$ has a finite orbit $F$ then the strong approximation conjecture cannot hold.  To see this consider more generally any finite orbit $F$ of the $\Gamma$ action on $\afine^3(\cx)$.  The coordinate of any $\xi$ in such an $F$ must lie in a cyclotomic field $L_n=\ratls(\zeta_n)$, where $\zeta_n$ is a primitive $n$-th root of $1$.  For if $\xi =(\xi_1, \xi_2, \xi_3)$ then $\ord( \rot(\xi_j))$ must be finite and hence 
\beq \label{p51} 3 \xi_j = t_j + t_j^{-1} \eeq with $t_j$ a root of unity.  If $n$ is the least common multiple of all the orders of all the $t_j$'s corresponding to the $\xi$'s in $F$ then $\xi_j \in L_n$ and hence $F \subset \afine^3(L_n)$.  In fact the $\xi_j$'s are all integral except possibly for denominators powers of $3$, so that $F \subset \afine^3 (\mathcal{O}_{L_n}^{(3)})$ where $\mathcal{O}_{L_n}^{(3)}$ is the ring of $S$-integers in $L_n$, with $S$ consisting of the primes dividing $3$.  If $p$ is a rational prime $(p \neq 3$) which splits completely in $L_n$ and $P$ is a prime of $L_n$ with $P |( p) $ then $\mathcal{O}_{L_n}^{(3)}/P \cong \fp (\cong \zed/p \zed)$.
The $\Gamma$ action of $\afine^3(\mathcal{O}_{L_n}^{(3)})$ factors through the reduction $\pi$ mod $P$
$$
\begin{CD}
\afine^3(\mathcal{O}_{L_n}^{(3)}) @> \Gamma >> \afine^3(\mathcal{O}_{L_n}^{(3)})\\
@VV \pi V      @VV\pi V \\
\afine^3(\mathcal{O}_{L}^{(3)}/P)@ >\bar{ \Gamma} >> \afine^3(\mathcal{O}_{L}^{(3)}/P)
\end{CD}
$$
and hence
\beq \label{p52}
 \bar{F}=\pi(F)\subset \afine^3(\fp), \,  \, \text{is} \, \, \bar{\Gamma}- \text{invariant}.
\eeq
Since $\Gamma$ preserves the level sets $X_k$:
\beq \label{p53}x_1^2+x_2^2+x_3^2-3 x_1 x_2x_3=k,
\eeq
any such $F$ is contained in $X_k(\mathcal{O}_{L}^{(3)})$ for a suitable $k$.  Thus for any such $F$, there is a positive density of $p$'s which split completely in $L_n$, and hence for which $\bar{F}\subset \mv_k(\fp)$ is a fixed size $\bar{\Gamma}$ -orbit ($|\bar{F}|\leq |F|$).
That is, the finite $\Gamma$-orbits in $\afine^3(\bar{\ratls})$ must be part of any description of the $\bar{\Gamma}$-orbits on $\afine^3(\fp)$, for $p$ large.

In our  setting of this paper, $k=0$ and we have (we thank E. Bombieri for this simple proof)
\begin{proposition} \label{prop51}
$\xx(\bar{\ratls})$ has no finite $\Gamma$-orbit.
\end{proposition}
\proof As in the discussion above, if $F$ is such an orbit and $\xi \in F$ then the $\xi_j$ satisfy \eqref{p51} with $t_j$ an $l_j$-th root of one.  The Markoff equation 
for $t_1, t_2, t_3$ becomes
\beq \label{p54}
(t_1+t_1^{-1})^2+(t_2+t_2^{-1})^2+(t_3+t_3^{-1})^2-(t_1+t_1^{-1})(t_2+t_2^{-1})(t_3+t_3^{-1})=0.
\eeq

Now \eqref{p54} has no solutions with $|t_j|=1$ (let alone being roots of unity) except for $t_j=\pm i$, $j=1,2,3$.  To see this note that if $$a=t_1+t_1^{-1} (=t_1+\bar{t_1}), \, \, b=t_2+t_2^{-1}, \, \, c=t_3+t_3^{-1}$$
then $a,b,c,$ lie in $[-2,2]$ and by the inequality of the geometric and arithmetic means 
\beq \label{p55}
0\leq a^2+b^2+c^2=|abc|\leq \frac{|a|^3+|b|^3+|c|^3}{3}\leq  \frac{2}{3}(a^2+b^2+c^2).\eeq
Hence the only solutions to \eqref{p54} correspond  to $a=b=c=0$ or $t_j=\pm i$.  In terms of $\xi_j$ this gives $\xi=(0,0,0)$, which is the only invariant set for the action of $\Gamma$ on $\mv(\bar{\ratls})$.
\qed

We remark that in the context of the general surfaces that are studied in Paper III, for example the surfaces $\mv_k$ in \eqref{p53} with $k\ne 0$, there can be a continuum of solutions to the analogue of equation \eqref{p54} with $|t_j|=1$.   However the solutions with $t_j$ a root of unity (with unspecified order) are still restricted to a finite number of nondegenerate solutions.  This follows from Lang's  $G_m$ conjecture,  see \cite{Lau83} and \cite{SA94} for proofs which give the solutions effectively.  In various special cases these finite $\bar{\ratls}$ orbits for the $\Gamma$-action correspond to the determination of the algebraic Painleve VI solutions (\cite{DM00}, \cite{LT14}); we leave the details to paper III .

Returning to the Markoff surface $\mv$, let $\xi=(\xi_1, \xi_2, \xi_3) \in \xx(p)$ with $\ord(\rot( \xi_j) )= l_j$ for $j=1,2,3$.  Let $n=\lcm(l_1,l_2,l_3)$ 
and $L_n=\ratls(\zeta_n)$ and let $\zeta_{l_1}, \zeta_{l_2}, \zeta_{l_3}$ be primitive roots of one respectively.
Let
\beq \label{p56}
\eta=(\zeta_{l_1}+\zeta_{l_2}^{-1})^2+(\zeta_{l_2}+\zeta_{l_2}^{-1})^2+(\zeta_{l_3}+\zeta_{l_3}^{-1})^2-
(\zeta_{l_1}+\zeta_{l_2}^{-1})(\zeta_{l_2}+\zeta_{l_2}^{-1})(\zeta_{l_3}+\zeta_{l_3}^{-1})\in \mathcal{O}_{L_n}.\eeq

According to Proposition \ref{prop51} unless
$l_1=l_2=l_3=2$  (i.e. $\zeta_{l_j}=\pm i$), $\eta \ne 0$.

Now $|\eta| \leq 20$ and hence 
\beq \label{p57}
\Norm(\eta)\leq 20^{\phi(n)}\leq 20^{n}.\eeq
If $P$ is a prime in $\mathcal{O}_{L_n}$ and $\eta \in P$, then $P | (\eta)$ and hence 
\beq \label{p58}
\N(P) \leq \Norm(\eta)\leq 20^n.\eeq
Put differently, if
$$\log_{20} \N(P) > n$$
then
\beq \label{p59}
\eta \ne 0 (\, \Mod\,  P). \eeq

For our point $\xi$ in $\xx(p)$, $3 \xi_j = \lambda_j+\lambda_j^{-1}$ with $\lambda_j$ in $\fp$ or $\fpp$ and $\lambda_j$ an $l_j$-th root of $1$, 
and $(l_1, l_2, l_3)\ne (2,2,2)$ since $\xi \ne (0,0,0)$.  If all the $\lambda_j$'s are in $\fp$ then $\ratls(\zeta_n)$ splits completely at $p$, that is there  is a 
prime $P$ dividing $(p)$ such that 
$$\mathcal{O}_{L_n}/P \cong \fp \, \, , \, \, \N(P)=p$$
and $\pi(\zeta_{l_j})=\lambda_j$ in $\mathcal{O}_{L_n}/P$ and $\eta \equiv 0 (\mod P)$.  Hence from \eqref{p59} we conclude that 
\beq \label{p510}
\log_{20}p \leq n.\eeq

If the field generated by $\lambda_j$'s (over $\fp$) is $\fpp$ then there is a prime $P$ of $\mathcal{O}_{L_n}$ dividing $(p)$ such that 
$$\mathcal{O}_{L_n}/P \cong \fpp \, \, , \, \, \N(P)=p^2$$
and $\pi(\zeta_{l_j})=\lambda_j$ in $\mathcal{O}_{L_n}/P$ and $\eta \equiv 0 (\, \Mod \,  P)$.  Hence again from \eqref{p59} we conclude that 
\beq \label{p511}
2\log_{20}p \leq n.\eeq
Hence in either case $n \geq \log_{20}p$ where $n=\lcm(l_1, l_2, l_3)$, and hence 
\beq \label{p512}
\max(l_1,l_2, l_3) \geq (\log_{20} p)^{\frac{1}{3}}.
\eeq
We have proven
\begin{proposition} \label{prop53}
Let $\xi \in \xx(p)$ have maximal order $l$, i.e. $\max(l_1, l_2, l_3)=l$ with $l_j = \ord(\rot (\xi_j))$, then $l\geq (\log_{20} p)^{\frac{1}{3}}$.
In particular,  any component $F$ of $\xx(p)$ satisfies
$$|F| \geq (\log_{20}p)^{\frac{1}{3}}.$$
\end{proposition}

\section{Proofs of Theorems  \ref{t1} and  \ref{t2}} \label{s6}

Proposition \ref{prop53} establishes the second part of Theorem \ref{t1} and combined with  the analysis  in sections \ref{s3} and \ref{s4}  yields a proof of the strong approximation conjecture if $p^2-1$ is not very smooth.  For example, the strong approximation conjecture is true for $\xx(p)$ if the prime $p$ satisfies
\beq \label{smooth}
\sum_{\substack{(\log p)^{\frac{1}{3}} \leq d \leq y\\d|(p^2-1)}} d^{\frac{2}{3}} < y; \, \, \text{for any} \, \, y.  \eeq

We do this  by using  the arguments and results in [Cha13] and [CKSZ14] concerning points $(x, y)$ on irreducible curves over $\mathbb F_p$ for which $\ord (x)+ \ord (y)$
is small (here $\ord (x)$ is the order of $x$ in $\mathbb F_p^*$).

\begin{theorem} \label{order}  Fix $d\in \zed_+$ and $\delta>0$.   There is an $\varepsilon>0$, $\varepsilon = \varepsilon(d, \delta)$, such that for  all primes $p \leq z$ ($z$  sufficiently  large)  with the exception of at most $z^{\delta}$ of them,  the following property holds.  Let  $f(x, y) \in \fp[x,y]$ be  of degree at most $d$ and not divisible by any polynomial of the form $\rho x^{\ga}y^{\gb} -1$ or $\rho y^{\gb}-x^{\ga}$ for any $\rho \in \bar{\fp}$ and integers $\ga$ and $\gb$.  Then all solutions $(x,y)\in (\bar{\fp} \times \bar{\fp})^{*}$ of $f(x,y) =0$ satisfy
\beq \label{e71}
\ord(x)+\ord(y) \geq p^{\varepsilon} \eeq except for at most $11 d^3 + d$ of them.\end{theorem}

\vspace{5mm}

\proof Theorem 1.2 in \cite{CKSZ14} establishes what we want except that $p^{\varepsilon}$ in \eqref{e71} is replaced by the stronger bound $p^{\ga(d)}$, with $$\ga(d)= \frac{2}{89 d^2 + 3d + 14},$$  while the exceptional set of primes is of zero density.  For our purpose the exponent in \eqref{e71} is allowed to be small and in exchange we want the exceptional set to be much smaller.  To this end we follow verbatum the discussion in Section 4  of \cite{CKSZ14}.  For $d$ fixed and $T$ a large parameter they show that there is a $U=U(d, T)$ which has at most  $O(T^{\frac{1}{\ga(d)}})$ prime factors   (their $\log T$ in the denominator is irrelevant for us) with the property: If $p$ does not divide $U$ and $f$ as in Theorem \ref{order} and $f(x,y)=0$ in $\fp$, then 
\beq \label{e72} \ord(x)+\ord(y) \geq T \eeq
except for at most $11  d^3+ d$ such $(x,y)$ in $\fp^{*} \times \fp^{*}$.    For our given $\delta>0$ and large parameter $z$ choose $T$  to be 
\beq \label{e73} T = z^{\delta \ga(d)}. \eeq
Then the number of primes $p$ with $p|U$ is $O(z^{\delta})$ , and if $p$ does not divide $U$ then Theorem \ref{order} holds with \eqref{e72} and \eqref{e73}, that is with $\varepsilon = \delta \ga(d)$.

\qed

\vspace{5mm}

To prove Theorem \ref{t2} we apply Theorem \ref{order} to the curves $f_{\sigma}(x,y)$ given by equation 
$$x+\frac{\sigma}{x}= y+\frac{1}{y}$$
with $\sigma \neq 1$.  If $(\log_{20} p)^{\frac{1}{3}} > 1000$, then according to Proposition \ref{prop53}, for any $\xi \in \xx(p)$ we have $\ord(\rot(\xi_{j_0}))$ is at least $1000$ for $j_0$ one of the $j$'s in $\{1, 2, 3\}$.  Hence if $p$ is not in the exceptional set in Theorem \ref{order} with $d=4$, then in  the typical equation  $x+\frac{\sigma}{x}= y+\frac{1}{y}$ corresponding to the orders of the rotations in the $\rot(\xi_{j_0})$ orbit, there is $(x, y)$ which is not one of the exceptional $11d^3 +d < 1000$ possible points.  For such $(x, y)$ the induced rotation has order at least $p^{\varepsilon}$ and hence $\xi$ is joined to the cage by the middle game.

\vspace{5mm}

Our methods fall short of dealing with all $p$, specifically for those rare $p$'s for which $p^2-1$ is very smooth.
The following hypothesis which is a strong variant of conjectures of M.C.~Chang and B.~Poonen [Cha13, Vol10] would suffice to deal with all large $p$'s.

\vspace{.1cm}

\noindent
{\bf Hypothesis:}
{\sl Given $d\in\mathbb N$, there is $\delta>0$ and $K=K(d)$ such that for $p$ large and $f(x, y)$ absolutely irreducible over $\mathbb F_p$ and of degree  $d$ and $f(x,
y)=0$ is not a translate of a subtorus of $(\bar{\mathbb F}_p^*)^2$, the set of $(x, y)\in (\mathbb F_p^*)^2$ for which $f(x, y)=0$ and $\max(\ord x, \ord y)\leq p^\delta$, is
at most $K$.}

\appendix
\section{}

Stepanov's auxiliary polynomial method \cite{Ste69} for bounding the number of solutons to equations like \eqref{p41} is quite flexible.
We demonstrate this for some special cases (the general case can be handled similarly).   The proposition below is an extension of the approach and bounds in \cite{HK00}
(where $S(x)=x$, $T(x)=1-x$ and $t_1=t_2$).

In what follows $S(x)$  and $T(x)$ are rational functions in $\fp(x)$ of total degree $d_1$ and $d_2$ respectively and with disjoint divisors; $e=d_1+d_2$ is fixed.

\begin{proposition} \label{a1} For $p$ a large prime, $t_1, t_2$ dividing $p-1$, $t_1\geq t_2$, let
$$Y=\{y\in \fp\, : \, S(y)^{t_1}=T(y)^{t_2}=1 \}.$$
Then if $t_1 \ll_{e} p^{1-\frac{1}{2e}}$,
$$|Y|\ll_{e} \min \{t_2, t_1 t_2^{-\frac{1}{4e}}\}.$$
\end{proposition}

\noindent \textbf{Remarks:}
\begin{enumerate}
\item The trivial bound is $O(t_2)$ so the Proposition gives an improvement (power saving) if $t_2 \geq t_1^{\frac{4e}{4e-1}}$.
\item If $h(\xi, \eta)=0$ is a plane curve of genus $0$ over $\fp$, then the Proposition gives an upper bound on the number of solutions with $\xi^{t_1} =\eta^{t_2} =1$
(cf. \cite{CZ13}).
\end{enumerate}

Applying Proposition \ref{a1} with $t_1=t_2$, $S(y)=y$, $T(y)=\frac{ay+b}{cy+d}$ yields
\begin{corollary} \label{a4} For $p$ large prime, $t|(p-1)$, $t\leq p^{\frac{3}{4}}$ and $U_t=\{y\in \fp \, :\, y^t=1\}$  the $t$-th roots of $1$,
$$|\sigma(U_t) \cap U_t|\ll t^{\frac{3}{4}}$$
for $\sigma \in \PGL_2(\fp)$, $\sigma \neq 1$.
\end{corollary}

\begin{corollary} \label{a5}For $t|(p-1)$, $t\leq p^{\frac{3}{4}}$, $b\in \fp$, $b\neq 1$,
$$|\{w, \rho \in \fp\, :\, w+w^{-1} = \rho + b \rho^{-1}, w^t=\rho^t=1\}| \ll t^{\frac{3}{4}}.$$
\end{corollary}
\proof Put $\rho w = \xi, \frac{w}{\rho}=\eta$, then $\xi^t=\eta^t=1$  and each such solution with $\xi=\frac{b \eta-1}{\eta-1}$ corresponds to at most two solutions $(w, \rho)$ above.   Applying Corollary \ref{a4} yields Corollary \ref{a5}. \qed

\bigskip

\noindent \textbf{Proof of Proposition \ref{a1}:}
First we need a generalization of Proposition 3.2 in \cite{VS12} where their common $t$ is replaced by $t_0, t_1, \dots, t_n$.  The result is the following Lemma, whose proof is the same
\begin{lemma} \label{a6} Let $t_0, t_1, \dots, t_n$, as well as $B$ and $J$ be integers, $p$ a large prime, and $\ga_1, \dots, \ga_n$ distinct elements in $\fp^{*}$.  Assume that 
$$\min(t_0, \dots, t_n) \geq \frac{1}{2}(n-1) B^{2n} +JB$$
and that 
$$p\geq (2nB+2) \max(t_0, t_1, \dots,  t_n).$$
Then
$$x^{a_i} x^{t_0 b_{0, i}} (x-\ga_1)^{t_1 b_{1,i}}\dots  (x-\ga_n)^{t_n b_{n,i}}$$
with $a_j \leq J$ and $b_{0,i}, \dots b_{n,i} \leq B$ are linearly independent in $\fp[x]$. \end{lemma}

Let $\ga_1, \dots, \ga_k \in \fp$ be distinct and $\nu_1, \nu_2, \dots \nu_k \in \zed$; set
$$R_{\nu}(x)=(x-\ga_1)^{\nu_1} \dots (x-\ga_k)^{\nu_k}.$$
For $m\geq 1$,
\begin{equation*} \begin{split} &\frac{d^m}{d x^m}\left[ R_{\nu}(x)\right]=\sum_{j_1+\dots j_k=m} \binom{m}{j} \frac{d^{j_1}}{dx^{j_1}}
\left[(x-\ga_1)^{\nu_1}\right] \dots   \frac{d^{j_k}}{dx^{j_k}}
\left[(x-\ga_1)^{\nu_k}\right] =\\
&\sum_{j_1+j_2 +\dots j_k =m} B_{m,j} (x-\ga_1)^{\nu_1-j_1} \dots (x-\ga_k)^{\nu_k-j_k}.\end{split} \end{equation*}
Hence
\beq \label{a1e}
\left[(x-\ga_1)(x-\ga_2) \dots(x-\ga_k\right]^m\frac{d^m}{d x^m} R_{\nu}(x)=R_{\nu}(x)P_{m, \nu}(x), \eeq
where $P_{m,\nu}$ is a polynomial of degree at most $km$.

Stepanov's polynomial method is based on constructing a polynomial which vanishes to high order on $Y$.
Let $\lambda_{a, b_1, b_2}$ be in $\fp$ with $0 \leq a \leq J$ and $0\leq b_j \leq B$ and form
\beq  \label{a2e}
\phi(x)=\sum \lambda_{a, b_1, b_2} x^{a} (S(x))^{t_1 b_1} (T(x))^{t_2 b_2}.\eeq
Write $S(x), T(x)$ in the form (we assume that both factor into linear factors $\fp[x]$):
\beq \label{a3e} 
\begin{aligned} &S(x)=\frac{A (x-\ga_1) \dots (x-\ga_t)}{(x-\gb_1)\dots (x-\gb_{\tau})}, \\ &T(x)=\frac{B (x-\gamma_1) \dots (x-\gamma_{\mu})}{(x-\delta_1) \dots (x-\delta_{\nu})}. \end{aligned} \eeq
For simplicity we assume that $S(x)$ and $T(x)$ are  square-free and we are assuming that the $\alpha, \beta, \gamma, \delta$'s are all distinct.   The constants $A$ and $B$ can be absorbed into the $\lambda$'s, so without loss of generality we can take $A=B=1$.
For $m\geq 0$
\beq \label{a4e}
\begin{split}
[(x-\ga_1) \dots (x-\ga_t)(x-\beta_1) \dots (x-\beta_{\tau}) \dots (x-\delta_{\nu})]^{m}& \frac{d^m}{d x^{m}}\left[ x^a  (S(x))^{t_1 b_1} (T(x))^{t_2 b_2}\right]=\\ &x^a S(x)^{t_1 b_1} T(x)^{t_2 b_2} P_{m,a, b_1, b_2}(x) \end{split} \eeq
with $P_m$ of degree at most $em$.
Hence for $x=y \in Y$ and $m\leq M$
\beq \label{a5e}\frac{d^m}{d x^m} \phi(x) |_{x=y} = \sum \lambda_{a, b_1, b_2} \frac{d^m}{d x^m}\left[ x^a S(x)^{t_1 b_1} T(x)^{t_2 b_2}\right]_{x=y}=
\sum_{a, b_1, b_2} \lambda_{a, b_1, b_2} y^a P_{m, a, b_1, b_2}(y),\eeq
by the definition of $Y$.  We can make \eqref{a5e} equal to $0$ for all $y$ in $Y$  by noting that $y^a P_m(y)$ is a polynomial of degree $J+em$. So \eqref{a5e} can be made $0$ with not all of the $\lambda_{a, b_1, b_2}$'s equal to $0$ and for all $m \leq M$ as long as 
\beq \label{a6e} (J+ eM) M <B^2 J. \eeq
Assuming that this is satisfied, we have $\phi(x)$ which is not identically zero and has degree (as a rational function) at most $J+eBt_1$.

Hence if $\phi(x)$ is not identically zero, then its order of vanishing on $Y$ is at least $M$ and hence
$$M |Y|\leq J + e B t_1$$  or \beq \label{a7e}|Y|\leq\frac{J+e Bt_1}{M}.\eeq
We now check that under suitable constraints on the sizes of parameters,  $\phi(x)$ does not vanish identically.
We have
$$\phi(x)= \sum \lambda_{a, b_1, b_2}x^a \frac{(x-\ga_1)^{t_1 b_1}  \dots (x-\ga_t)^{t_1 b_1}}{(x-\beta_1)^{t_1b_1} \dots (x-\beta_{\tau})^{t_1b_1}}
\frac{(x-\gamma_1)^{t_2 b_2}  \dots (x-\gamma_{\mu})^{t_2 b_2}}{(x-\delta_1)^{t_2b_2} \dots (x-\delta_{\nu})^{t_2b_2}}, $$
consequently
\begin{equation*}
\begin{split} &(x-\beta_1)^B \dots (x-\beta_{\tau})^B (x-\delta_1)^B \dots (x-\delta_{\nu})^B \phi(x) =\\
&\sum \lambda_{a, b_1, b_2}x^a (x-\ga_1)^{t_1 b_1}  \dots (x-\ga_t)^{t_1 b_1}(x-\beta_1)^{(B-b_1)t_1} \dots (x-\beta_{\tau})^{(B-b_1)t_1}\cdot\\
&\cdot (x-\gamma_1)^{t_2 b_2}  \dots (x-\gamma_{\mu})^{t_2 b_2}(x-\delta_1)^{(B-b_2) t_2} \dots (x-\delta_{\nu})^{(B-b_2) t_2}. 
\end{split}
\end{equation*}
Now the monomials appearing in the last expression are linearly independent over $\fp[x]$ according to Lemma \ref{a6} as long as 
\beq \label{a8e} 
\begin{aligned} &p \geq (2 e B+2) t_1 \\ & t_2 \geq \frac{1}{2} e B^{2 e} + J B^{e}. 
\end{aligned}\eeq
Thus,  if \eqref{a8e} and \eqref{a6e} hold, so does \eqref{a7e}.

Choose $J \leq M$ and $$M^2=C_e B^2 J.$$
Then $$M=\sqrt{C_e} B \sqrt{J}$$ and $M\geq J$ iff $J \ll_{e} B^2$.
Now choose $B+t_1^{\frac{1}{2e}}$ and if $t_1 \ll p^{1-\frac{1}{2e}}$ then \eqref{a8e} and \eqref{a6e} hold and $$|Y| \ll_{e}t_1 t_2 ^{-\frac{1}{4e}}.$$
This completes the proof of Proposition \ref{a1}.

\noindent
{\bf Acknowledgements:} 
While working on this paper, the authors were supported, in part, by the following NSF DMS  awards: 1301619 (Bourgain), 1603715  (Gamburd), 1302952 (Sarnak).


\begin{thebibliography}{BHWVGRS11}

\bibitem[BS02]{BS02} F. Beukers and C. J. Smyth, Cyclotomic points on curves, Number theory for the millenium (Urbana, Illinois, 2000), I, A.K. Peters, (2002), 67-85.



\bibitem[Bom07]{Bom07} E.~Bombieri, Continued fractions and the Markoff tree, Expo. Math. 25 (2007), no 3,
187--213.

\bibitem[Bou12] {Bou12} J.~Bourgain, A modular Szemeredi-Trotter theorem for hyperbolas, C.R. Acad. Sci. Paris Ser 1,
350 (2012), 793--796.




\bibitem[BG08]{BG08} J.~Bourgain and A.~Gamburd, Uniform expansion bounds for Cayley graphs of $\SL_2(\mathbb
F_p)$, Ann. Math. 167 (2008), 625--642.





\bibitem[BGS10]{BGS10} J.~Bourgain, A.~Gamburd and P.~Sarnak, Affine linear sieve, expanders and sum product,
Invent. Math. 179 (2010), 559--644.







\bibitem[BGS16]{BGS16} J. Bourgain, A. Gamburd,  P. Sarnak,  Markoff triples and strong approximation,  C. R. Math. Acad. Sci. Paris 354 (2016), no. 2, 131-135.



 






\bibitem[Cas78]{Cas78} J. W. S. Cassels, \emph{Rational Quadratic Forms},
Academic Press, 1978.






\bibitem[Cha13]{Cha13}Chang, Mei-Chu Elements of large order in prime finite fields. Bull. Aust. Math. Soc. 88 (2013) 169-176.


\bibitem[CKSZ14]{CKSZ14} M-C.~Chang, B.~Kerr, I.~Shparlinski and U.~Zannier,
Elements of large orders on varieties over prime finite fields,  J. Theor.  Nombres Bordeaux 26 (2014) 579-594.





\bibitem[CS13]{CZ13} P. Corvaja,  U. Zannier,  Greatest common divisors of $u-1, v-1$ in positive characteristic and rational points on curves over finite fields,  J. Eur. Math. Soc. (JEMS) 15 (2013) 1927-1942.



\bibitem[DKS13]{DKS13} C. D'Andrea, T. Krick and M. Sombra, Heights of varieties in multiprojective spaces and arithmetic Nullstellensatz, Annales Sci. de l'ENS, 46, (2013),549-627.



\bibitem[DM00]{DM00}B. Dubrovin, M.  Mazzocco,  Monodromy of certain Painleve-VI transcendents and reflection groups, Invent. Math. 141 (2000) 55-147.






\bibitem[Fro13]{Fro13} G.~Frobenius, \"Uber die Markoffschen Zahlen, Akad. Wiss. Berlin, (1913), 458--487.




\bibitem[Gol03]{Gol03} W.~Goldman, The modular group action on real
$SL(2)$-characters of a one-holed torus, Geom. and Top. Vol. 7 (2003),
443--486.







\bibitem[HK00]{HK00} R.~Heath-Brown and S.~Konyagin, New bounds for Gauss sums derived from $k$-th powers and
for Heilbronn's exponential sum, QUAR. J MATH (2000), 221--235.











\bibitem[Lau83]{Lau83} M.~Laurent, Exponential diophantine equations, CR Acad Sc. 296 (1983), 945--947.


\bibitem[Mar79]{Mar79} A. Markoff,  Sur les formes quadratiques binaires ind\'{e}finies, Math. Ann. 15 (1879) 381--409.



\bibitem[Mar80]{Mar80} A. Markoff,  Sur les formes quadratiques binaires ind\'{e}finies, Math. Ann. 17 (1880) 379--399.

\bibitem[LT14]{LT14} O. Lisovyy and Y. Tykhyy,  Algebraic solutions of the sixth Painleve  equation, Journal of Geometry and Physics, \textbf{85} (2014) 124-163.









\bibitem[MW13]{MW13} D.~McCullough and M.~Wanderley,  Nielsen equivalence of generating pairs in $SL(2, q)$,
Glasgow Math. J. 55 (2013), 481--509.






\bibitem[SA94]{SA94}P. Sarnak and S. Adams,  Betti numbers of congruence groups, with an appendix by Z. Rudnick, Israel J. Math, 88, 1994, 31-72.


\bibitem[Sch04] {Sch04}Schmidt, W. M \emph{Equations over finite fields: an Elementary approach}, Kendrick Press 2004.




\bibitem[Ste69]{Ste69} S.A.~Stepanov, The number of points of a hyperelliptic curve over a prime field,
MATH USSR-IZV 3:5 (1969), 1103--1114.




\bibitem[Suz82]{Suz82}
M. Suzuki, \emph{Group Theory I}, Springer-Verlag,
Berlin-Heidelberg-New York, 1982.






\bibitem[Vol10]{Vol10}Voloch,  Elements of high order on finite fields from elliptic curves, Bull. Aust. Math. Soc. 81 (2010), 425-429.


\bibitem[VS12]{VS12}  Vyugin, I. V., Shkredov, I. D., Sb. Math. 203(2012), 844-863.


\bibitem[Wei41]{Wei41}A.~Weil, On the Riemann Hypothesis in function fields,  Proc. Nat. Acad.  Sci. USA 27 (1941),
345--347.


 
 




\end{thebibliography}
\end{document}